\documentclass[12pt,a4j]{article}

\usepackage{amsmath,amssymb}
\usepackage{amsthm}
\usepackage[nobysame]{amsrefs}
\usepackage{latexsym}
\usepackage[all]{xy}
\usepackage{mathrsfs}
\usepackage[dvipdfmx]{graphicx}

\newtheorem{Thm}{Theorem}
\newtheorem{Prop}[Thm]{Proposition}
\newtheorem{Lem}[Thm]{Lemma}
\newtheorem{Cor}[Thm]{Corollary}
\newtheorem{Def}[Thm]{Definition}
\theoremstyle{definition} % The sentences in the following environments will not become italic automatically by this command. 

\newtheorem{Rem}[Thm]{Remark}

\newcommand{\Proof}{{\it Proof.}\quad}
\newcommand{\QED}{{\unskip\nobreak\hfil\penalty50\quad\null\nobreak\hfil
$\square$\parfillskip0pt\finalhyphendemerits0\par\medskip}}

\newfont{\eufm}{eufm10}
\newfont{\eufms}{eufm8}
\newfont{\msam}{msam10}
\newfont{\msams}{msam8}
\newfont{\msbm}{msbm10}

\renewcommand{\geq}{\mathrel{\hbox{\msam\symbol{"3D}}}}

\newcommand{\be}{\begin{equation}}
\newcommand{\ee}{\end{equation}}

\newcommand{\Z}{{\mathbf Z}}
\newcommand{\N}{{\mathbf N}}
\newcommand{\cl}{\mathrel{:}}

\newcommand{\ra}{\rightarrow}

\newcommand{\calS}{\mathcal{S}}
\newcommand{\calC}{\mathcal{C}}
\newcommand{\calP}{\mathcal{P}}

\makeatletter
\def\mapstofill@{%
   \arrowfill@{\mapstochar\relbar}\relbar\rightarrow}
\newcommand*\xmapsto[2][]{%
   \ext@arrow 0395\mapstofill@{#1}{#2}}
\makeatother

\title{Preprojective roots and cycle graphs}
\date{}
\author{Yuji Komatsu}

\begin{document}
\maketitle

\begin{abstract}
%Let $(W,S)$ be a Coxeter system and $\Phi$ be the root system of $(W,S)$.
%It is known that every root is either positive or negative.
%In particular, we consider the Coxeter group whose Coxeter graph is cyclic. 
We study $c$-preprojective roots for a Coxeter element $c$ of infinite Coxeter group $W$. 
In particular, we consider the case when any positive root is $c$-preprojective
 for some Coxeter element $c$. 
In this paper, by assuming that the Coxeter graph of $W$ is cyclic,
 we establish that any positive root is $c$-preprojective
 for some Coxeter element $c$ if and only if $W$ is an affine Coxeter group. 
\end{abstract}

Crystallographic Coxeter groups are used in several recent papers
 for representations of quivers or more general finite dimensional algebras.
 %see for example []
On the other hand, Pelley and Kleiner used representations of quivers to prove that
 the powers of a Coxeter element are reduced
 in an infinite irreducible crystallographic Coxeter group \cite{Pelley},
 and then Speyer proved the result for a general Coxeter group \cite{Speyer},
 using combinatorics of \cite{Pelley} and stripping out the quiver theory. 
These papers show that certain results on representations of quivers have analogs
 in the theory of Coxeter groups. 
The way from quivers to Coxeter groups goes through the notion of root in view of the results of Kac \cite{Kac}. 
A root of a Coxeter group is an analog of a real root of a quiver,
 and each positive real root of a quiver is the dimension vector
 of a unique up to isomorphism indecomposable representation. 
Therefore we view a positive root of a Coxeter group $W$ as an analog of an indecomposable representation. 
% and view a finite set of positive roots as an analog of a representation
% that need not be indecomposable but has no isomorphic direct summands.

Among the indecomposable representations of an acyclic quiver that correspond to real roots,
 the most important are preprojective and preinjective representations introduced
 by Bernstein, Gelfand, and Ponomarev \cite{BGP} as those annihilated by a power of the Coxeter functor. 
The analog of the Coxeter functor is a Coxeter element,
 so Igusa and Schiffler \cite{I-S} fix a Coxeter element $c \in W$
 and define a $c$-preprojective root as a positive root sent to a negative root
 by a positive power of $c$.  
The element $c$ determines a unique acyclic orientation of the Coxeter graph $\Gamma$ of $W$
 and thus turns it into a quiver. 
The inverse Coxeter element, $c^{-1}$, yields the opposite quiver,
 so the $c^{-1}$-preprojective roots are analogs of preinjective representations.

Kleiner shows that the Coxeter group is finite if and only if,
 for all Coxeter element $c$, each positive root is $c$-preprojective \cite{Klein}. 
His result implies that when the Coxeter group $W$ is infinite,
 for each Coxeter element $c$ there exists a positive root $\alpha$
 which is not $c$-preprojective. 
Here we may ask whether there is another Coxeter element $c'$
 such that such positive root $\alpha$ is $c'$-preprojective. 
%Let  $W$ be an infinite Coxeter group. 
%Then his result says that, for each Coxeter element $c \in W$,
% there exists a positive but not $c$-preprojective root $\alpha$. 
%However, $\alpha$ may be $c'$-preprojective for another Coxeter element $c' \in W$. 
Indeed, if $W$ is an affine Cxeter group of $\widetilde{A}_{n}$ type
 then, for any positive root $\alpha$, there exists a Coxeter element $c \in W$
 such that $\alpha$ is $c$-preprojective root
 (see Proposition \ref{affine to cup}). 
The Coxeter graph of $\widetilde{A}_{n}$ type is cyclic. 
For a Coxeter group with a cyclic graph, we obtain a criterion for the 
group to be affine in terms of the $c$-preprojecivity of positive roots.
%When Coxeter graph is cycle, we obtain criterion for classification of Coxeter groups
% according to the property of $c$-preprojective roots: 

\begin{Thm}
\label{affine preproj.}
If the Coxeter graph of $W$ is cyclic, then the following statements are equivalent: 
\begin{description}
  \item[($1$)] $W$ is an affine Coxeter group, 
  \item[($2$)] each positive root is $c$-preprojective for some Coxeter element $c$.
\end{description}
\end{Thm}

In Sections $1$ and $2$, we recall definitions and results about Coxeter groups. 
The notion of depth of positive roots is important here. 
We show the necessary condition for a Coxter group to be infinite and non-affine
 in Corollary \ref{growth for dp}, where the depth is used. 
In Sections $3$ and $4$, we suppose that a Coxeter graph is cyclic. 
We define Coxeter elements denoted by $c_i^k$
 and focus on $c_i^k$-preprojective roots. 
Their properties are used in the proof of Theorem \ref{affine preproj.}. 
The statement of Theorem \ref{affine preproj.} is separated into
 Proposition \ref{cup to affine} and Proposition \ref{affine to cup}.

%Theorem $1$ is proved by Proposition \ref{cup to affine} in Section $3$
% and Proposition \ref{affine to cup} in Section $4$.
%Finally, we conjecture characterization that 
% the Coxeter group $W$ is affine if and only if each positive root is $c$-preprojective
% for some Coxeter element $c$.

We will report that the same statement of Proposition \ref{affine to cup} remains valid
 for affine Coxeter groups of different types in a next opportunity.

%%%%%%%%%%  %%%%%%%%%%

\section{ROOT SYSTEM OF A COXETER SYSTEM}
%1    intro

We recall some facts about Coxeter groups,
 based on \cite{Klein, BB}. 
Several properties which we use in this paper are stated as propositions. 
% to prove Theorem \ref{affine preproj.}

Let $(W,S)$ be a {\it Coxeter system},
 where $W$ is a group generated by $S =\{ s_1,s_2,\ldots,s_n \}$. 
The defining relations are $(s_is_j)^{m_{i,j}}=1$ for $s_i$, $s_j\in S$ if $m_{i,j}<\infty$,
 where the $m_{i,j}$ are the entries of a {\it Coxeter matrix}
 $(m_{i,j}=m(s_i,s_j))_{i,j=1,2,\ldots,n}$. 
The Coxeter matrix is a symmetric $n\times n$ matrix with $m_{i,j}\in \Z \cup \{\infty\}$;
 $m_{i,i}=1$ for all $s_i\in S$; and $m_{i,j}>1$ whenever $i\neq j$. 
Denote by $\Gamma=(\Gamma_0,\Gamma_1,m)$ the {\it Coxeter graph} of $(W,S)$
 with the set of vertices $\Gamma_0=S$, the set of edges $\Gamma_1$,
 and the map $m\cl\Gamma_1\ra\{3,4,\ldots,\infty\}$. 
There exists one, and only one,
 edge joining vertices $s_i$ and $s_j$ if and only if $2<m_{i,j}$. 
The map $m$ assigns $m_{i,j}$ to the edge joining vertices $s_i$ and $s_j$. 
We call $(W,S)$ {\it irreducible} if its Coxeter graph is connected. 
Each element $w\in W$ can be written as $w=t_1t_2\cdots t_k$, where $t_i \in S$. 
If $k$ is minimal among all such expressions for $w$,
 then $k$ is called the {\it length} of $w$ (written as $\ell(w)=k$)
 and the word $t_1t_2\cdots t_k$ is called a {\it reduced expression} for $w$.

Let $V$ be an $n$-dimensional real vector space, 
 with the basis $\Pi=\{ \alpha_i \}_{i=1,2,\ldots,n}$ parametrized by $S$ so that
 $\alpha_i$ corresponds to $s_i$. 
%Also, we write as $\alpha_i=\alpha_{s_i}$. %since $\mathrm{card}(S)=n=\mathrm{card}(\Phi)$. 
Let $(\cdot \mid \cdot)$ be the symmetric bilinear form such that
 $ (\alpha_i \mid \alpha_j)\cl=-\cos\pi/m_{i,j}$, where $\pi/\infty=0$. 
Then $W$ acts on $V$ by
\[  s_i(v) \cl= v-2(\alpha_i \mid v)\alpha_i \ \text{for $v \in V$}. \]
A vector $\beta=w(\alpha_i)$ for some $w\in W$ and $\alpha_i\in\Pi$ is called a {\it root} ;
 a root $\beta$ is {\it positive} (respectively, {\it negative})
 if all of its coordinates are nonnegative (respectively, nonpositive)
 when expressed in the basis $\Pi$. 
It is well known that every root is either positive or negative. 
Let $\Phi$ be the {\it root system} of $(W,S)$, which is the set of all roots,
 and let $\Phi^+$ (respectively, $\Phi^-$) be the set of all positive roots (respectively, negative roots). 
Note that $\Phi^+=-\Phi^-$.

%%%%%%%%%% Prop %%%%%%%%%%
\begin{Prop}\cite[Proposition\ 4.4.6]{BB}
\label{l to root}
For $w\in W$, let $w=t_1t_2\cdots t_k$ be a reduced expression. 
Then 
\begin{align*}
 \Phi^+\cap w^{-1}\bigl(\Phi^-\bigr)&=\{ \beta\in\Phi^+ \mid w(\beta)\in\Phi^- \}\\
                &=\{ \alpha_{t_k}, t_k(\alpha_{t_{k-1}}),\ldots, t_kt_{k-1}\cdots t_2(\alpha_{t_1}) \}.
\end{align*}
Hence $\ell(w)=\mathrm{card}(\Phi^+\cap w^{-1}\bigl(\Phi^-\bigr))$. 
\end{Prop}

The {\it depth} of $\beta \in \Phi^+$ is the minimal $k$
 such that $\beta = t_1t_2\cdots t_{k-1}(\alpha_i)$ for some $t_j\in S$ and $\alpha_i\in \Pi$;
 written as $\mathrm{dp}(\beta)=k$. 
%We have how the depth changes when acting on positive roots by a generator : 

%%%%%%%%%% Prop %%%%%%%%%%
\begin{Prop}\cite[Lemma\ 4.6.2]{BB}
\label{depth change}
For $s_i\in S$ and $\beta \in \Phi^+\setminus\{\alpha_i\}$, $s_i(\beta)$ is also positive root
 and has the following equality: 
\[
\mathrm{dp}(s_i(\beta))=\begin{cases} 
         \mathrm{dp}(\beta)-1  &  \text{if $(\alpha_i \mid \beta) > 0$},\\
         \mathrm{dp}(\beta)     &  \text{if $(\alpha_i \mid \beta) = 0$},\\
         \mathrm{dp}(\beta)+1  &  \text{if $(\alpha_i \mid \beta) < 0$}. 
         \end{cases}
\]
\end{Prop}

Let $c$ be an element of $W$,
 given by $c=s_{\sigma(1)}s_{\sigma(2)}\cdots s_{\sigma(n)}$
 for some permutation $\sigma \in \calS_n$. 
This $c$ is called a {\it Coxeter element }. 
It is known that the length of Coxeter element is $n$. 
In particular, if $c= t_1t_2 \cdots t_n$ is a reduced expression for a Coxeter element $c$
 then there exists a permutation $\sigma \in \calS_n$ such that $t_j = s_{\sigma(j)}$,
 $j=1, 2, \ldots, n$. 
Denote by $\calC$ the set of all Coxeter elements. 
A Coxeter element $c$ gives an orientation to $\Gamma$
 so that every arrow $s_{\sigma(i)}\ra s_{\sigma(j)}$ satisfies $i<j$. 
Then $S$ become a poset (partially ordered set) by setting $s\le_c t$
 if there exists a path from $s$ to $t$. 
We denote the poset by $(S,c)$. 
Note that this correspondence of a Coxeter element to a poset on $S$ is injective.

We fix a Coxeter element $c\in W$. 
A positive root $\beta$ is called {\it $c$-preprojective}
 if $c^{\mu}(\beta)$ becomes a negative root for some positive integer $\mu$.
Denote by $\calP(c)$ the subset of $\Phi^+$ consisting of all $c$-preprojective roots.

%%%%%%%%%% Prop %%%%%%%%%%
\begin{Prop}\cite[Theorem\ 1]{Speyer}
\label{c^n}
Let W be an infinite, irreducible Coxeter group. 
Let $c=t_1t_2\cdots t_{n}$ be a Coxeter element, where $t_i=s_{\sigma(i)}\in S$. 
Then, for any positive integer $\mu$, the expression
 $c^{\mu}=(t_1\cdots t_{n})(t_1\cdots t_{n})\cdots (t_1\cdots t_{n})$
 ($\mu$ times) is a reduced expression for $c^{\mu}$. 
\end{Prop}

%%%%%%%%%% Cor %%%%%%%%%%
\begin{Cor}\cite{Klein}
\label{Per(c)}
In the setting of Proposition \ref{c^n},
\[ \Phi^+\cap c^{-\mu}\bigl(\Phi^-\bigr) =
 \Bigl(\Phi^+\cap c^{-\mu+1}\bigl(\Phi^-\bigr)\Bigr) \cup c^{-1}\Bigl(\Phi^+\cap c^{-\mu+1}\bigl(\Phi^-\bigr)\Bigr) \]
for any positive integer $\mu$. 
This implies 
\[ \calP(c) = \bigcup_{\mu = 0}^{\infty} c^{-\mu}\Bigl(\Phi^+\cap c^{-1}\bigl(\Phi^-\bigr)\Bigr) .\]
\end{Cor}

%\begin{Rem}
%\label{orbit of Per(c)}
%For any element $w \in W$, 
%\[ w\Bigl( \Phi^+\cap w^{-1}\bigl(\Phi^-\bigr) \Bigr) = w\bigl(\Phi^+\bigr) \cap \Phi^- 
%                                   = - \Bigl( (w^{-1})^{-1}\bigl(\Phi^-\bigr) \cap \Phi^+ \Bigr) . \]
%Hence if $c$ is a Coxeter element then
% $c^\mu(\beta)$ belongs to $\calP(c) \cup -\calP(c^{-1})$
% for any $\mu \in \Z$, $\beta \in \calP(c)$. 
%\end{Rem}

%%%%%%%%%%  %%%%%%%%%%

\section{THE CLASSIFICATION OF A COXETER SYSTEM}
%2    criterion

In the following, we assume the Coxeter system is irreducible. 
We can classify Coxeter systems into three types as finite, {\it affine}, and others,
 according to the signature of the bilinear form $(\cdot \mid \cdot)$: 
\begin{description}
  \item[($1$)] if $(\cdot \mid \cdot)$ is positive definite,
                           then $(W,S)$ is finite, 
  \item[($2$)] if $(\cdot \mid \cdot)$ is positive semidefinite
                           with $0$ a simple eigenvalue,
                           then $(W,S)$ is affine, 
  \item[($3$)] none of the above conditions hold. 
\end{description}
%The above characteristics have some equivalences. 
Indeed, $(W,S)$ is a finite Coxeter system if and only if $W$ is a finite group;
 if $(W,S)$ is a affine Coxeter system then $W$ is isomorphic to some affine reflection group.

We know another criterion for this classification introduced by Terragni \cite{Terra}. 
This criterion depends on the ({\it exponential}) {\it growth rate}. 
For a non-negative integer $r$, let $W_r\cl=\{ w \in W \mid \ell(w)=r \}$. 
Then the growth rate of $(W,S)$ is defined as
 $\omega(W,S) \cl= \limsup_{r\ra\infty} \sqrt[\leftroot{-1} \uproot{2} r] {\mathrm{card}(W_r)}$.

%%%%%%%%%% Thm %%%%%%%%%%
\begin{Thm}\cite{Terra}
\label{growth}
Let $(W,S)$ be a Coxeter system. 
Then: 
\begin{description}
  \item[($1$)] $(W,S)$ is a finite Coxeter system if and only if $\omega(W,S)=0$, 
  \item[($2$)] $(W,S)$ is an affine Coxeter system if and only if $\omega(W,S)=1$, 
  \item[($3$)] $(W,S)$ is a Coxeter system of other type if and only if $\omega(W,S)>1$. 
\end{description}
\end{Thm}

In this paper, we use the property of growth of depth in place of Theorem \ref{growth}. 
Let $\Phi^+(r) \cl= \{ \beta \in \Phi^+ \mid \mathrm{dp}(\beta)=r \}$ for a positive integer $r$.

%%%%%%%%%% Cor %%%%%%%%%%
\begin{Cor}
\label{growth for dp}
%Let $(W,S)$ be a Coxeter system. 
If $\omega(W,S)>1$ then $\limsup_{r\ra\infty} \mathrm{card}(\Phi^+(r)) = \infty$. 
\end{Cor}
\Proof
Assume that $\limsup_{r\ra\infty} \mathrm{card}(\Phi^+(r)) < \infty$. 
Then there exists a positive integer $N$ such that $\mathrm{card}(\Phi^+(r)) < N$
 for any positive integer $r$. 
We consider the following injective map:
\[
M_r \cl W_r \ni w \mapsto (w(\alpha_1) , w(\alpha_2) , \ldots , w(\alpha_n)) \in \Phi^n .
\]
If $w(\alpha_i)$ is positive root then $\mathrm{dp}(w(\alpha_i)) \le r+1$
 since $w$ has a reduced expression $w = t_1t_2\cdots t_r$ where $t_j \in S$. 
In the case of $w(\alpha_i) \in \Phi^-$,
 we obtain $\alpha_i \in \Phi^+ \cap w^{-1}\bigl( \Phi^- \bigr)$ and
\[ ws_i(\alpha_i) = w(-\alpha_i) = -w(\alpha_i) \in \Phi^+ . \]
Since if $\beta \in \Phi^+\setminus \{\alpha_i\}$ then $s_i(\beta)$ is also positive root,
\begin{align*}
 \Phi^+ \cap (ws_i)^{-1}\bigl( \Phi^- \bigr) &= \Phi^+ \cap s_iw^{-1}\bigl( \Phi^- \bigr) \\
		&= s_i\Bigl( \bigl(\Phi^+\setminus \{\alpha_i\}\bigr) \cup \{-\alpha_i\} \Bigr) \cap s_iw^{-1}\bigl( \Phi^- \bigr) \\
		&= s_i\Bigl( (\bigl(\Phi^+\setminus \{\alpha_i\}\bigr) \cup \{-\alpha_i\}) \cap w^{-1}\bigl( \Phi^- \bigr) \Bigr) \\
		&= s_i\Bigl( \bigl(\Phi^+ \cap w^{-1}\bigl( \Phi^- \bigr) \bigr)\setminus \{\alpha_i\} \Bigr) . 
\end{align*}
Thus $\ell(ws_i) = \mathrm{card}(\bigl(\Phi^+ \cap w^{-1}\bigl( \Phi^- \bigr) \bigr)\setminus \{\alpha_i\}) = r-1$
 by Proposition \ref{l to root}, this implies $\mathrm{dp}(-w(\alpha_i)) = \mathrm{dp}(ws_i(\alpha_i)) \le r$.

Hence the image of $M_r$ is included in
 $\bigl( (\bigcup_{k=1}^r -\Phi^+(k) )\cup(\bigcup_{k=1}^{r+1} \Phi^+(k) ) \bigr)^n$. 
Therefore we infer $\mathrm{card}(W_r) \le N^n(2r+1)^n$ for any positive integer $r$, that is,  $\omega(W,S)\le1$. 
It is contradiction. 
\QED

%%%%%%%%%%  %%%%%%%%%%

\section{CYCLIC COXETER GRAPHS}
%3      maintheorem

We consider the case that the Coxeter graph of $(W,S)$ is cyclic with $n= \mathrm{card}(S)\geq3$. 
Suppose that $m_{i,j}$ is equal to $2$ if and only if $j=i+1$ or $j=i-1$ ($n+1=1$ $0=n$). 
Note that $m_{i,j}=2$ implies the following properties:
\begin{itemize}
  \item $s_i$ and $s_j$ commute, 
  \item $s_i(\alpha_j)=\alpha_j$. 
\end{itemize}
For simplicity, we denote $m_j=m_{j,j+1}=m_{j+1,j}$. 
The Coxeter graph of $(W,S)$ is illustrated in Figure 1. 
\begin{figure}[h]
  \centering
  \includegraphics[width=4cm]{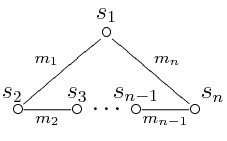}
  \caption{The cyclic Coxeter graph}
\end{figure}

Then $W$ is an infinite group. 
In particular, $(W,S)$ is affine if and only if $m_j=3$ for $j=1, 2, \ldots, n$. 
(This affine Coxeter system is called of $\widetilde{A}_{n-1}$ type. )

%%%%%%%%%% Def %%%%%%%%%%
\begin{Def}
\label{c^i_k}
For $\mu \in \Z$, there exists a unique expression $\mu=\lambda n+\nu$
 such that $\lambda , \nu \in \Z$ and $1\le \nu \le n$. 
We denote $[\mu]\cl=\nu$. 

For $i$, $k \in \Z$ with $i<k<i+n$,
 we define a Coxeter element $c_{[i]}^{[k]}$ as
\be
c_{[i]}^{[k]}\cl=s_{[i]}(s_{[i+1]}s_{[i+2]}\cdots s_{[k-1]})(s_{[i+n-1]}s_{[i+n-2]}\cdots s_{[k+1]})s_{[k]} . 
\ee
Since this definition depends only on the value of $[i]$ and $[k]$,
 we can write $c_i^k=s_{i}(s_{[i+1]}\cdots s_{[k-1]})(s_{[i-1]}\cdots s_{[k+1]})s_{k}$ for any two integers $1 \le i$, $k \le n$. 
\end{Def}

It is clear that $(c_i^k)^{-1}=c^i_k$. 
The poset $(S , c_i^k)$ (defined after Proposition \ref{depth change})
 has the greatest element $s_k$, and the least element $s_i$. 
Note that, in the case of cyclic Coxeter graphs,
 the number of maximal elements is $1$ if and only if
 the number of minimal elements is $1$. 
Thus, for any $c\in \calC$, if the poset $(S , c)$ has a greatest element
 then $c$ is one of the $c_i^k$.

\begin{Prop}
\label{conjugate}
For any Coxeter element $c$, there exist a element $w$ of $W$
 and some integers $i$, $k$, $1 \le i , k \le n$ such that $wcw^{-1}=c_i^k$. 
\end{Prop}

This proposition follows from the next Lemma. 
Let $M(S,c)$ be the set of all maximal elements of $(S,c)$
 and let $(s)_c\cl=\{ t \in S \mid t \le_c s \}$
 for a Coxeter element $c \in \calC$ and a element $s \in S$.

\begin{Lem}
\label{maxmin}
Let $c\in W$ be a Coxeter element with $\mathrm{card}(M(S,c)) \ge 2$. 
For any maximal element $s \in M(S,c)$, there exists an element $w \in W$
 satisfying the following three properties: 
\begin{description}
  \item[($1$)] $wcw^{-1}$ is also a Coxeter element, 
  \item[($2$)] $M(S,wcw^{-1}) \subset M(S,c)$, 
  \item[($3$)] $s \notin M(S,wcw^{-1})$ or $\mathrm{card}((s)_{wcw^{-1}}) < \mathrm{card}((s)_c)$ holds. 
\end{description}
\end{Lem}
\Proof
By the definition of $(S,c)$, there exist integers $\lambda \le \mu \le \nu$ such that 
\[ (s)_c=\{ s_{[\lambda]}, s_{[\lambda+1]}, \ldots , s_{[\mu]}=s, \ldots , s_{[\nu]} \}, \]
\[ s_{[\lambda]} \le_c s_{[\lambda+1]} \le_c \ldots \le_c s \ge_c s_{[\mu+1]} \ge_c \ldots \ge_c s_{[\nu]}. \] 
Then $s_{[\lambda]}$ and $s_{[\nu]}$ are minimal elements in $(S,c)$. 
These facts implies the reduced expression (\ref{sp.exp}) for $c$. 
Let $c=t_1t_2\cdots t_n$ be a reduced expression. 
Now if $t_j=s_{[\lambda]}$ (resp. $s_{[\nu]}$) then $s_{[\lambda-1]}, s_{[\lambda+1]}$
 (resp. $s_{[\nu-1]}, s_{[\nu+1]}$) are not elements of $\{ t_1, t_2, \ldots, t_{j-1} \}$
 since $s_{[\lambda]}$ (resp. $s_{[\nu]}$) is a minimal element. 
Thus we can take a reduced expression for $c$:
 \be
\label{exp1}
c = s_{[\lambda]}s_{[\nu]}t_1t_2 \cdots t_{j-1}t_{j+1} \cdots t_{h-1}t_{h+1} \cdots t_n
\ee
for some integers $1\le j < h \le n$. 
Similarly, if $t_j=s_{[\mu]}$ then $s_{[\mu-1]}, s_{[\mu+1]}$ are not elements of
 $\{ t_{j+1}, t_{j+2,} \ldots, t_n \}$ since $s_{[\mu]}$ is a maximal element. 
%Similarly, if $s_h \le_c t_j$ for some $s_h \in S$, $j=1, 2, \ldots, n$
% then $s_h \notin \{ t_{j+1}, t_{j+2,} \ldots, t_n \}$. 
By the method similar to above for $s_{[\mu]}$ and the reduced expression (\ref{exp1}),
 we obtain a reduced expression:
\[ c = s_{[\lambda]}s_{[\nu]}t'_1t'_2 \cdots t'_{n-3}s_{[\mu]} . \]
similarly to the above, we can write as
 \be
\label{sp.exp}
c=s_{[\lambda]}s_{[\nu]}x(s_{[\lambda+1]}s_{[\lambda+2]}\cdots s_{[\mu-1]})(s_{[\nu-1]}s_{[\nu-2]}\cdots s_{[\mu+1]})s_{[\mu]}
\ee
where $x = s_{\sigma([\nu+1])}s_{\sigma([\nu+2])} \cdots s_{\sigma([\lambda-1+n])}$
 for some permutation $\sigma$ of $S \setminus (s)_c$. 
%with $\ell(x)=n-\nu+\lambda-1$. 
Let $w=(s_{[\lambda+1]}s_{[\lambda+2]}\cdots s_{[\mu-1]})(s_{[\nu-1]}s_{[\nu-2]}\cdots s_{[\mu+1]})s_{[\mu]}$. 
Then 
\[ wcw^{-1} = (s_{[\lambda+1]}s_{[\lambda+2]}\cdots s_{[\mu-1]})(s_{[\nu-1]}s_{[\nu-2]}\cdots s_{[\mu+1]})s_{[\mu]}s_{[\lambda]}s_{[\nu]}x \]
is also a Coxeter element.

Let $c'$ be this Coxeter element $wcw^{-1}$. 
%We denote this Coxeter element $wcw^{-1}$ by $c'$.  
We shall prove that $c'$ has properties $(2)$, $(3)$ of Lemma \ref{maxmin}
 in three different cases:
 when $\mu=\nu-1$, $\mu=\lambda+1$, and others. 
If $\mu=\nu-1$ then %we obtain
\[ c' = s_{[\lambda+1]}s_{[\lambda+2]}\cdots s_{[\nu-1]}s_{[\nu]}s_{[\lambda]}x .\]
%in the some method as above. 
Thus ordering of $(s)_c$ in $(S,c')$ is
\[ s_{[\lambda]} \ge_{c'} s_{[\lambda+1]} \le_{c'} s_{[\lambda+2]} \le_{c'} \ldots
 \le_{c'} s_{[\nu-1]} = s \le_{c'} s_{[\nu]} .\]
Note that the ordering of $S\setminus (s)_c$ depend on a reduced expression for $x$.  
Therefore $M(S,c')= M(S,c)\setminus \{ s \}$. 
Similarly, if $\mu=\lambda+1$ then $M(S,c')= M(S,c)\setminus \{ s \}$.

In other cases, we obtain a reduced expression:
\[ c' = s_{[\lambda+1]}s_{[\nu-1]}s_{[\lambda]}s_{[\nu]}x(s_{[\lambda+2]}\cdots s_{[\mu-1]})(s_{[\nu-2]}\cdots s_{[\mu+1]})s_{[\mu]} \]
in the some method as above. 
Now the ordering of $(s)_c$ in $(S,c')$ is
\[ s_{[\lambda]} \ge_{c'} s_{[\lambda+1]} \le_{c'} \ldots \le_{c'} s \ge_{c'} \ldots \ge_{c'} s_{[\nu-1]} \le_{c'} s_{[\nu]} .\]
Thus $(s)_{c'}=\{ s_{[\lambda+1]}, s_{[\lambda+2]}, \ldots , s_{[\nu-1]} \}$ and $M(S,c')= M(S,c)$. 
\QED

%Note that the ordering of $S\setminus (s)_c$ depend on a reduced expression for $x$.  
%By the equation, we obtain that
% the ordering of $S\setminus\{  s_{[\lambda+1]}, \ldots , s_{[\nu-1]} \}$ in $(S,c)$
% is equal to the ordering of $S\setminus\{  s_{[\lambda+1]}, \ldots , s_{[\nu-1]} \}$ in $(S,c')$. 

{\it Proof of Proposition \ref{conjugate}.}\quad
Fix a maximal element $s_k$ of $(S,c)$. 
By using Lemma \ref{maxmin} some times, we obtain a sequence $\{ w_j \}_{1\le j \le h}$
 and a Coxeter element $c'=w_h\cdots w_1cw_1^{-1}\cdots w_h^{-1}$
 such that $s_k$ is a greatest element in $(S,c')$. 
% and the three property of Lemma \ref{maxmin} hold
% for each Coxeter element $w_{j-1}\cdots w_1cw_1^{-1}\cdots w_{j-1}^{-1}$. 
Then there exists a least element $s_i$  in $(S,c')$. 
This means that $(S,c')$ is equal to $(S,c_i^k)$, that is, $c'=c_i^k$. 
\QED

For any Coxeter element $c \in W$, Corollary \ref{Per(c)} says that
 if $\beta$ is a $c$-preprojective root then $c^{-1}(\beta)$ is also a $c$-preprojective root. 
Moreover, if $c$ is equal to $c_i^k$ for some integers $1 \le i , k \le n$
 then $\beta$ and $c^{-1}(\beta)$ have the next relation:

%%%%%%%%%% Prop %%%%%%%%%%
\begin{Prop}
\label{preproj.}
Let $i$ and $k$ be two integers with $1 \le i , k \le n$. 
Then $\beta=\sum_{j=1}^n b_{j}\alpha_{j} \in \calP(c_i^k)$ has following properties:
\begin{description}
  \item[(1)] $b_h \le b_j$ when $s_h \le_{c_i^k} s_j$, 
  \item[(2)] $\mathrm{dp}(c^i_k(\beta)) > \mathrm{dp}(\beta)$
\end{description}
\end{Prop}
\Proof
We list each element of $\Phi^+\cap c^i_k\bigl(\Phi^-\bigr)$
 since $\calP(c_i^k) = \bigcup_{\mu = 0}^{\infty} (c^i_k)^\mu (\Phi^+\cap c^i_k\bigl(\Phi^-\bigr))$ (Corollary \ref{Per(c)}). 
By Proposition \ref{l to root} and the definition of $c_i^k$, 
\begin{align*}
\Phi^+\cap c^i_k\bigl(\Phi^-\bigr) =
       \{ \alpha_k , &s_k(\alpha_{[k+1]}) , s_ks_{[k+1]}(\alpha_{[k+2]}) , s_ks_{[k+1]} \cdots s_{[i-2]}(\alpha_{[i-1]}) , \\ 
          &s_k(\alpha_{[k-1]}) , s_ks_{[k-1]}(\alpha_{[k-2]}) , s_ks_{[k-1]} \cdots s_{[i+2]}(\alpha_{[i+1]}) , \\
          &s_k(_{[k+1]} \cdots s_{[i-1]})(s_{[k-1]} \cdots s_{[i+1]})(\alpha_i)\} .
\end{align*}
Let $(x_{1,h} , x_{2,h} , \dots , x_{n,h})$ be the representation of
 $s_ks_{[k+1]} \cdots s_h(\alpha_{[h+1]}) \in \Phi^+\cap c^i_k\bigl(\Phi^-\bigr)$
 with respect to $\Pi=\{ \alpha_1, \ldots, \alpha_n \}$. 
Then, for an integer $g$ with $s_h \le_{c_i^k} s_g \le_{c_i^k} s_k$, 
\[ s_gs_{[g+1]} \cdots s_h(\alpha_{h+1}) = \sum_{s_h \le_{c_i^k} s_j \le_{c_i^k} s_g} x_{j,h}\alpha_j \]
by the definition of group action. 
Thus 
\[ x_{j,h} = \begin{cases} 
         1                                      & \text{ if $j=h+1$ }, \\
         2\cos\frac{\pi}{m_j}x_{[j+1],h} & \text{ if $ s_h \le_{c_i^k} s_j \le_{c_i^k} s_g$ }, \\
         0                                      &  \text{others}.                  
         \end{cases} \]  %cfrac
Since $m_j \ge 3$ for all $j$, a $c_i^k$-preprojective root $s_ks_{[k+1]} \cdots s_h(\alpha_{[h+1]})$ has the property $(1)$. 
Similarly, we can prove that each element of $\Phi^+\cap c^i_k\bigl(\Phi^-\bigr)$
 has the property $(1)$.

%For $s_k(_{[k+1]} \cdots s_{[i-1]})(s_{[k-1]} \cdots s_{[i+1]})(\alpha_i)$,
% let $(y_1 , x_2 , \dots , y_n})$ be the representation of this with respect of $\Pi$. 

Note that $\calP(c_i^k) = \{ (c_k^i)^\mu(\gamma) \mid \mu \in \N \cup \{0\} ,
 \gamma \in \Phi^+\cap c^i_k\bigl(\Phi^-\bigr) \}$. 
We prove inductively Proposition \ref{preproj.} by next Lemma.

%%%%%%%%%% Lem %%%%%%%%%%
\begin{Lem}
\label{depth}
%Let $\lambda , \nu$ be integrs with $\lambda \le \nu<\lambda+n$. 
Let $\beta=\sum_{j=1}^n y_{j}\alpha_{j} \in V$
 with %$ I \cl= \{ j \mid y_j \neq 0 \}= \{ [\lambda],[\lambda+1],\ldots,[\nu]\}$,
 $y_j \ge 0$ for all $j$. 
Suppose $y_h \le y_j$ for all $s_h , s_j \in S$ satisfying $s_h \le_{c_i^k} s_j$. 
Then $c_k^i(\beta)=\sum_{j=1}^{n} z_{j}\alpha_{j} $ has the following properties: 
\begin{description}
%  \item[(1)] $\displaystyle J \cl= \{ j \mid z_j \neq 0 \}$\\
%  \ \ \       $= \begin{cases} 
%          \{ [\lambda-1],[\lambda],\ldots,[\nu+1]\}  &  \text{if} \ [i-1], [i+1] \notin I , \\
%         \{ 1, 2, \ldots , n\}    &    \text{others},  
%          \end{cases}$ 
  \item[(1)] $z_h \le z_j$ when $s_h \le_{c_i^k} s_j$, 
  \item[(2)] $z_j \ge y_{[j-1]} , y_j , y_{[j+1]} $ for $j = 1, 2, \ldots, n$. 
\end{description}
\end{Lem}
By using inductively Lemma \ref{depth},
 we obtain that each $c_i^k$-preprojective root has the property $(1)$ of the Proposition \ref{preproj.}.

For a $c_i^k$-preprojective root $\beta=\sum_{j=1}^n b_{j}\alpha_{j}$
 and integers $h = 1, 2, \ldots , n$, let $c_k^i(\beta)=\sum_{j=1}^n b'_{j}\alpha_{j}$ and
\[ \beta_h = \begin{cases} 
		s_{[h-1]} \cdots s_{[i+2]}s_{[i+1]}s_i(\beta) & \text{if $s_i \le_{c_i^k} s_h \le_{c_i^k} s_{[k-1]}$ }, \\
		s_{[h+1]} \cdots s_{[i-2]}s_{[i-1]}(s_{[k-1]} \cdots s_{[i+2]}s_{[i+1]})s_i(\beta) & \text{if $s_{[i-1]} \le_{c_i^k} s_h \le_{c_i^k} s_k$ }.
          \end{cases} \] 
Note that $\beta_i = \beta$, $s_k(\beta_k) = c_k^i(\beta)$ and
\be 
\label{s_hb_h} 
s_h(\beta_h) = \begin{cases} 
		\beta_{[h+1]} & \text{if $s_h \le_{c_i^k} s_{[k-1]}$ and $s_h \neq s_{[k-1]}$ }, \\
		\beta_{[i-1]} & \text{if $h=[k-1]$ }, \\
		\beta_{[h-1]} & \text{if $s_{[i-1]} \le_{c_i^k} s_h \le_{c_i^k} s_{[k+1]}$ }.
          \end{cases} 
\ee
Then, by the definition of group action, 
\[ \beta_h = \begin{cases} \displaystyle 
		\sum_{\substack{ s_j \le_{c_i^k} s_h \\ s_j \neq s_h }} b'_j\alpha_j + b_h\alpha_h + \sum_{s_j \not\le_{c_i^k} s_h} b_j\alpha_j & \text{if $s_i \le_{c_i^k} s_h \le_{c_i^k} s_{[k-1]}$ }, \\ \displaystyle 
		\sum_{s_j \le_{c_i^k} s_{[k-1]}} b'_j\alpha_j +\sum_{\substack{ s_j \le_{c_i^k} s_h \\ s_j \neq s_h }}b'_j\alpha_j + \sum_{s_h \le_{c_i^k} s_j} b_j \alpha_j & \text{if $s_{[i-1]} \le_{c_i^k} s_h \le_{c_i^k} s_k$ }.
          \end{cases} \] 
By the difference between each representation of $\beta_h$ and $s_h(\beta_h)$
 with respect to $\Pi$, we obtain $2(\alpha_h \mid \beta_h) = b_h - b'_h$
 for all $h = 1, 2, \ldots , n$. 
Lemma \ref{depth} says that
\[ 2(\alpha_h \mid \beta_h) \le 0 \ \text{for all $h = 1, 2, \ldots , n$}. \]
By Proposition \ref{depth change}, this inequality means
\be
\label{inequ of depth}
\mathrm{dp}(s_h(\beta_h)) \ge \mathrm{dp}(\beta_h) \ \text{for all $h = 1, 2, \ldots , n$}.
\ee
In particular, $\mathrm{dp}(s_h(\beta_h)) = \mathrm{dp}(\beta_h)$ if and only if
 the positive root $s_h(\beta_h)$ is equal to the positive root $\beta_h$. 
%\begin{align*}
%     & \mathrm{dp}(c_k^i(\beta)) = \mathrm{dp}(s_k(\beta_k)) \\
%\ge & \mathrm{dp}(\beta_k) = \mathrm{dp}(s_{[k+1]}(\beta_{[k+1]})) \\
%     & \vdots \\
%\ge & \mathrm{dp}(\beta_{[i-2]}) = \mathrm{dp}(s_{[i-1]}(\beta_{[i-1]})) \\
%\ge & \mathrm{dp}(\beta_{[i-1]}) = \mathrm{dp}(s_{[k-1]}(\beta_{[k-1]})) \\
%\ge & \mathrm{dp}(\beta_{[k-1]}) = \mathrm{dp}(s_{[k-2]}(\beta_{[k-2]})) \\
%     & \vdots \\
%\ge & \mathrm{dp}(\beta_{[i+1]}) = \mathrm{dp}(s_i(\beta_i)) \\
%\ge & \mathrm{dp}(\beta_i) = \mathrm{dp}(\beta) .
%\end{align*}

The relations (\ref{s_hb_h}) and (\ref{inequ of depth}) imply that
\[ \mathrm{dp}(\beta) = \mathrm{dp}(\beta_i) \le \mathrm{dp}(s_i(\beta_i)) = \mathrm{dp}(\beta_{[i+1]}) \le
\ldots \le \mathrm{dp}(s_k(\beta_k)) = \mathrm{dp}(c_k^i(\beta)). \]
Assume $\mathrm{dp}(c_k^i(\beta)) = \mathrm{dp}(\beta)$ for some $\beta \in \calP(c_i^k)$.
Then $ \mathrm{dp}(s_h(\beta_h)) = \mathrm{dp}(\beta_h) $ for $h = 1, 2, \ldots , n$. 
Therefore we obtain $(c_i^k)^{-1}(\beta) = c_k^i(\beta) = \beta$,
 it is contradictory to $\beta \in \calP(c_i^k)$. 
Hence the property $(2)$ of Proposition \ref{preproj.} holds for any $c_i^k$-preprojective root. 
\QED

%\[ \beta_j = s_{[j-1]} \cdots s_{[i+2]}s_{[i+1]}s_i(\beta) \]  
%if $j \in \{ i , [i+1] , \ldots , [i+p]=[k-1]\}$ , $0\le p < n$, 
%\[ \beta_j =	s_{[j+1]} \cdots s_{[i-2]}s_{[i-1]}(s_{[k-1]} \cdots s_{[i+2]}s_{[i+1]})s_i(\beta) \]
%if $j \in \{ k , [k+1] , \ldots , [k+q]=[i-1]\}$ , $0\le q < n$. 

{\it Proof of Lemma \ref{depth}.}\quad
By the definition of group action,
\[ z_j=\begin{cases} 
		2\cos\frac{\pi}{m_{[i-1]}} y_{[i-1]} - y_i + 2\cos\frac{\pi}{m_i} y_{[i+1]}
 &  \text{if $j=i$ }, \\
		2\cos\frac{\pi}{m_{[j-1]}} z_{[j-1]} - y_j + 2\cos\frac{\pi}{m_j} y_{[j+1]}
 &  \text{if $s_{[i+1]} \le_{c_i^k} s_j \le_{c_i^k} s_{[k-1]}$ }, \\
		2\cos\frac{\pi}{m_{[j-1]}} y_{[j-1]} - y_j + 2\cos\frac{\pi}{m_j} z_{[j+1]}
 &  \text{if $s_{[i-1]} \le_{c_i^k} s_j \le_{c_i^k} s_{[k+1]}$ }, \\
		2\cos\frac{\pi}{m_{[k-1]}} z_{[k-1]} - y_k + 2\cos\frac{\pi}{m_k} z_{[k+1]}
 &  \text{if $j=k$ }. \\
          \end{cases} \]
Note that $2\cos\pi/m_h \ge 1$ for $h=1, 2, \ldots, n$. 
When $j=i$, 
\[ z_i = 2\cos\frac{\pi}{m_{[i-1]}} y_{[i-1]} + ( 2\cos\frac{\pi}{m_i} y_{[i+1]} - y_i )
 \ge y_{[i-1]} \ge y_i \]
by the supposition. 
Similarly, we obtain 
\[ z_i = ( 2\cos\frac{\pi}{m_{[i-1]}} y_{[i-1]} - y_i ) + 2\cos\frac{\pi}{m_i} y_{[i+1]} \ge y_{[i+1]} . \]

When $s_{[i+1]} \le_{c_i^k} s_j \le_{c_i^k} s_{[k-1]}$,
 we assume $z_{[j-1]} \ge y_{[j-2]}, y_{[j-1]}, y_j$.
Then 
\[ z_j \ge z_{[j-1]} \ge y_{[j-1]}, y_j \ \text{and} \ z_j \ge y_{[j+1]} \]
 in the some method as above. 
By induction, we can prove 
\[ z_j \ge z_{[j-1]}, y_{[j-1]}, y_j, y_{[j+1]} \ \text{for} \ 
s_{[i+1]} \le_{c_i^k} s_j \le_{c_i^k} s_{[k-1]} . \] 
Similarly, we obtain 
\[ z_j \ge z_{[j+1]}, y_{[j-1]}, y_j, y_{[j+1]} \ \text{for} \ 
s_{[i-1]} \le_{c_i^k} s_j \le_{c_i^k} s_{[k+1]} . \]

Therefore it is clear that 
\[ z_k \ge z_{[k-1]} \ge y_{[k-1]}, y_k , \]
\[ z_k \ge z_{[k+1]} \ge y_{[k+1]} . \]
\QED

%We will prove that the condition $(2)$ of Theorem \ref{affine preproj.} is sufficient for
% the condition $(1)$ of Theorem \ref{affine preproj.}

Let $\Phi^+\cap c^i_k\bigl(\Phi^-\bigr) = \{ \gamma_1, \gamma_2, \ldots, \gamma_n \}$
 and recall $\Phi^+(r) = \{ \beta \in \Phi^+ \mid \mathrm{dp}(\beta)=r \}$. 
Proposition \ref{preproj.} says that
 $\mathrm{dp}((c_k^i)^\mu(\gamma_j)) = \mathrm{dp}((c_k^i)^\nu(\gamma_j)) $ if and only if
 $\mu$ is equal to $\nu$,
  where $1 \le j \le n$, $\mu$ and $\nu$ are nonnegative integers . 
Thus 
\be
\label{<n}
\mathrm{card}(\calP(c_i^k) \cap \Phi^+(r)) \le n \ \text{for any} \ r \in \N
\ee
 by $\calP(c_i^k) = \bigcup_{j=1}^n \{ \gamma_j, c_k^i(\gamma_j), (c_k^i)^2(\gamma_j), \ldots \}$.

%%%%%%%%%% Prop %%%%%%%%%%
\begin{Prop}
\label{cup to affine}
Let $(W,S)$ be a Coxeter system such that Coxeter graph is cyclic. 
If $\Phi^+=\bigcup_{c\in\calC}\calP(c)$ then $(W,S)$ is affine. 
\end{Prop}
\Proof
Assume that $\Phi^+=\bigcup_{c\in\calC}\calP(c)$ and $(W,S)$ is a Coxeter system of other type. 
We shall implies a contradiction. 
Note that $\calC$ is finite set by the definition of Coxeter element. 
By Corollary \ref{growth for dp}, for any numbers $N, L$, there exist a Coxeter element $c$
 and a positive integer $r$ such that $\mathrm{card}(\calP(c) \cap \Phi^+(r)) > N$ and $r >L$.

For $c\in\calC$, we take numbers $L_c$ as following:
\[ L_c = \min \{ \ell(w) \mid w \in W, wcw^{-1} = c_i^k \ \text{for some integers} \ 1 \le i, k \le n \} . \]
Let $L_0 = \max \{L_c\}_{c\in\calC}$. 
Then there exist a Coxeter element $c_0$ and a positive integer $r_0$ such that
\be
\label{L_0}
 \mathrm{card}(\calP(c_0) \cap \Phi^+(r_0)) >2n(L_0+1) \ \text{and} \ r_0 > L_0 . 
\ee 
For this $c_0$, there exist a element $w_0$ of $W$ and some integers $1 \le i_0 , k_0 \le n$
 such that $w_0c_0w_0^{-1}=c^{k_0}_{i_0}$ and $\ell(w_0) \le L_0 < r$.

Let $\beta\in \calP(c_0) \cap \Phi^+(r_0)$. 
Then there exists a positive integer $\lambda$ such that $c_0^\lambda(\beta) \in \Phi^-$,
 that is, $(c_0^{-1})^\lambda$ sends a positive root $-c_0^\lambda(\beta)$
 to a negative root $-\beta$. 
By Proposition \ref{preproj.} for $-c_0^\lambda(\beta) \in \calP(c_0^{-1})$,
 there exists a positive integer $\mu_0 \ge \lambda$ such that
% c_0^{\mu_0}(\beta) \in \Phi^- \ \text{and} \
\[ \mathrm{dp}(-c_0^{\mu_0}(\beta)) > \ell(w_0) . \] 
Therefore $w_0(\beta), w_0(-c_0^{\mu_0}(\beta)) \in \Phi^+$
 by Proposition \ref{depth change}. 
(Note that the depth of positive root $\alpha$ is equal to $1$ if and only if
 $\alpha$ is one of $\alpha_1, \ldots, \alpha_n$. )
This implies that $w_0(\beta)$ is $c^{k_0}_{i_0}$-preprjective root since
\[ (c^{k_0}_{i_0})^{\mu_0}w_0(\beta) = (w_0c_0w_0^{-1})^{\mu_0}w_0(\beta) 
                                                     = w_0c_0^{\mu_0}(\beta) . \]
In particular, by Proposition \ref{depth change},
\[ r_0 - L_0 \le r_0 - \ell(w_0) \le \mathrm{dp}(w_0(\beta))
 \le r_0 + \ell(w_0) \le r_0 + L_0 . \]

Hence we obtain the injective map
\[ w_0 \cl \calP(c_0) \cap \Phi^+(r_0) \ra
 \bigcup_{r=r_0-L_0}^{r_0+L_0} \Bigl(\calP(c^{k_0}_{i_0}) \cap \Phi^+(r)\Bigr) . \] 
On the other hand, by the inequalities (\ref{<n}) and (\ref{L_0}),
\begin{align*}
2n(L_0+1) &< \mathrm{card}(\calP(c_0) \cap \Phi^+(r_0)) \\
 &\le \mathrm{card} (\bigcup_{r=r_0-L_0}^{r_0+L_0} \Bigl(\calP(c^{k_0}_{i_0}) \cap \Phi^+(r)\Bigr))\\
           %= \sum_{r=r_0-L_0}^{r_0+L_0}\Bigl(\calP(c^{k_0}_{i_0}) \cap \Phi^+(r)\Bigr) \\
 &\le n(2L_0+1) 
\end{align*}
 which is a contradiction. 
\QED

%%%%%%%%%%  %%%%%%%%%%

\section{ROOT SYSTEM OF $\widetilde{A}_{n-1}$ TYPE}
%maintheorem

Let $n\ge3$ and $(W,S)$ be a Coxeter system of $\widetilde{A}_{n-1}$ type. 
Then
\[
s_i(\alpha_j)=\begin{cases} 
         -\alpha_j              &  \text{if $\lvert i-j \rvert = 0$},\\
         \alpha_j+\alpha_i    &  \text{if $\lvert i-j \rvert = 1, n-1$},\\
         \alpha_j                &  \text{if $2 \le \lvert i-j \rvert \le n- 2$}. 
         \end{cases}
\]
Thus if $\mathrm{dp}(s_i(\alpha_j)) > \mathrm{dp}(\alpha_j)$ for $i \neq j$ then
 $s_i(\alpha_j) = \alpha_k + \alpha_{[k+1]}$ for some $k=1, 2, \ldots, n$.

Let $\beta=\sum_{j=\lambda}^{\nu}\alpha_{[j]}$
 for some integers $\lambda \le \nu \le \lambda+n-1$. 
Then
\[
s_i(\beta)=\begin{cases} 
         -\beta = -\alpha_i             &  \text{if $i = [\lambda] = [\nu]$},\\
         \beta - \alpha_i    &  \text{if $[\lambda] \ne [\nu], [\nu+1]$ and $ i = [\lambda], [\nu] $},\\
         \beta+\alpha_i & \text{if $[\lambda-1] \neq [\nu], [\nu+1]$ and $i=[\lambda-1], [\nu+1]$},\\
         \beta + 2\alpha_i               &  \text{if $i=[\lambda-1] = [\nu+1]$}, \\
         \beta    &  \text{others. }
         \end{cases}
\]
In particular, if $\beta=\sum_{j=1}^{n}\alpha_{j}$ then $s_i(\beta) = \beta$ for any $i$,
 that is, $\sum_{j=1}^{n}\alpha_{j}$ is not a root for $(W,S)$. 
In the case of $i=[\lambda-1] = [\nu+1]$,
 we can represent $\beta = \sum_{j=1}^{n}\alpha_{j} - \alpha_i$
 and $s_i(\beta) = \sum_{j=1}^{n}\alpha_{j} + \alpha_i$. 
Therefore, for $i=1, 2, \ldots, n$,
 if $\beta \neq \alpha_i, \sum_{j=1}^{n}\alpha_{j} - \alpha_i$
 then $s_i(\beta)$ is also represented by $\sum_{j=\lambda'}^{\nu'}\alpha_{[j]}$
 for some integers $\lambda' \le \nu' \le \lambda'+n-1$.

For $\mu \in \Z$,
 we denote by $(\mu)_n$ the vector $\sum_{i=1}^n \mu\alpha_i \in V$. 
From the above,
 $\sum_{j=\lambda}^{\nu}\alpha_{[j]}$ and $(1)_n + \alpha_i$ are positive roots
 for any integers $\lambda \le \nu < \lambda+n-1$ and $i=1, 2, \ldots, n$. 
It is clear that these positive roots are represented by $w(\alpha_k)$
 for some $w \in W$, $k=1, 2, \ldots, n$ by the definition of root. 
Hence $(1)_n + \sum_{j=\lambda}^{\nu}\alpha_{[j]}$ and $(2)_n + \alpha_i$ are also positive roots
 since $w((1)_n+\alpha_k)=(1)_n+w(\alpha_k)$. 
By the induction, we obtain the next Lemma. 

 %$\Phi=\{ (\mu)_n+\sum_{i=\lambda}^{\nu}\alpha_{[i]} \mid
% \mu, \lambda, \nu \in \Z, \lambda \le \nu <\lambda+n-1 \}$,

%%%%%%%%%% Lem %%%%%%%%%%
\begin{Lem}
\label{affine's root}
If $(W,S)$ be a Coxeter system of $\widetilde{A}_{n-1}$ type
 then $\Phi^+=\{ (\mu)_n+\sum_{j=\lambda}^{\nu}\alpha_{[j]} \mid
 \mu, \lambda, \nu \in \Z_{\ge 0}, \lambda \le \nu <\lambda+n-1 \}$. 
\end{Lem}

%%%%%%%%%% Prop %%%%%%%%%%
\begin{Prop}
\label{affine to cup}
Let $(W,S)$ is a Coxeter system of $\widetilde{A}_{n-1}$ type. 
Then $\Phi^+$ is covered by $\{ \calP(c) \}_{c\in\calC}$. 
\end{Prop}
\Proof
We consider the $c_{k+1}^k$-preprojective roots for $1\le k \le n$ ($n+1=1$). 
First, by Proposition \ref{l to root}, we obtain 
\[ \Phi^+\cap c^{k+1}_k\bigl(\Phi^-\bigr) = \{ \alpha_k, \alpha_{[k-1]}+\alpha_k, \ldots , \sum_{j=0}^{n-2}\alpha_{[k-j]}, (1)_n+\alpha_{k} \}. \]
By simple calculation, for $0 \le \lambda \le n-2$,
\[c^{k+1}_k(\sum_{j=0}^{\lambda}\alpha_{[k-j]}) =
\begin{cases} \displaystyle 
		(1)_n+\sum_{j=0}^{\lambda+1}\alpha_{[k-j]}    &  \text{if $0 \le \lambda \le n-3$},\\
		(2)_n+\alpha_k   &  \text{if $\lambda = n-2$}.
         \end{cases} \]
Therefore
\begin{align*}
 \calP(c_{k+1}^k) &= \bigcup_{\nu=0}^\infty (c^{k+1}_k)^\nu(\{ \alpha_k, \alpha_{[k-1]}+\alpha_k, \ldots , \sum_{j=k+2-n}^{k}\alpha_{[j]}, (1)_n+\alpha_{k} \}) \\
                       &= \{ (\mu)_n+\sum_{j=0}^{\lambda}\alpha_{[k-j]} \mid
 \mu, \lambda \in \Z_{\ge0} \ \text{and} \  \lambda \le n-2 \}.
\end{align*}

Hence we obtain $\Phi^+= \bigsqcup_{k=1}^n \calP(c_k^{k+1})=\bigcup_{c\in\calC}\calP(c)$ by Lemma \ref{affine's root}. 
\QED

%Proposition \ref{cup to affine}, \ref{affine covering} implies the main theorem. 

%\begin{Thm}[Main]
%\label{affine preproj.}
%Let $\calC$ be the set of Coxeter elements of $W$. 
%If the Coxeter graph of $W$ is cyclic, then the following statements are equivalent. 
%\begin{description}
%  \item[($1$)] $W$ is an affine Coxeter group, 
%  \item[($2$)] $\Phi^+=\bigcup_{c\in\calC}\calP(c)$. 
%\end{description}
%\end{Thm}

%%%%%%%%% Rem %%%%%%%%%%
\begin{Rem}
\label{remark}
Even if $n=2$, the statement of Proposition \ref{affine to cup} holds. 
Indeed a Coxeter system of $\widetilde{A}_1$ type is determined by $S=\{s_1,s_2\}$ and $m_{1,2}=m_{2,1}=\infty$. 
Then we easily obtain the following result: 
\[
s_i(\alpha_j)=\begin{cases} 
         -\alpha_i   &  \text{if $i = j$},\\
         \alpha_j + 2\alpha_i   &  \text{if $i \neq j$},
         \end{cases}
\]
\[ \Phi^+=\{ \alpha_1+(\mu)_2, \alpha_2+(\mu)_2 \mid \mu \in \Z_{\ge 0} \}, \]
\[ \calP(s_1s_2)= \{ \alpha_2+(\mu)_2 \mid \mu \in \Z_{\ge 0} \}, \]
\[ \calP(s_2s_1)= \{ \alpha_1+(\mu)_2 \mid \mu \in \Z_{\ge 0} \}. \]
% Proposition \ref{affine to cup} is deduced. 
\end{Rem}

\end{document}